\documentclass{amsart}

\usepackage{amsmath,amsfonts,amssymb,amsthm,amscd,latexsym,euscript,xypic}
\usepackage[all]{xy}

\newtheorem {theorem1}{Theorem}[section]
\newtheorem {theorem}[theorem1]{Theorem}

\newtheorem {lemma}[theorem1]{Lemma}
\theoremstyle{definition}

\theoremstyle{remark}

\newcommand{\calD}{\ensuremath{\mathcal{D}}}
\newcommand{\calS}{\ensuremath{\mathcal{S}}}
\newcommand{\calX}{\ensuremath{\mathcal{X}}}

\newcommand{\cat}[1]{\ensuremath{\EuScript #1}}

\DeclareMathOperator{\map}{\textup{map}}
\DeclareMathOperator{\Map}{\ensuremath{\textup{Map}}}

\newcommand{\colim}{\ensuremath{\mathop{\textup{colim}}}}

\newcommand{\rarrow}{\rightarrow}
\newcommand{\Id}{\ensuremath{\textup{Id}}}

\newcommand{\Vopenkas}{Vop\v enka's }

\newcommand{\HocatM}{\ensuremath{\textup{Ho}{\cat M}}}
\newcommand{\identity}{\ensuremath{\textup{id}}}
\newcommand{\ArrM}{\ensuremath{\textup{Arr}{\cat M}}}
\newcommand{\Sing}{\ensuremath{\textup{Sing}\,}}
\newcommand{\hbot}{\ensuremath{{\mathrm{h}\bot}}}

\begin{document}

\SelectTips{cm}{10}

\title[The orthogonal subcategory problem in homotopy theory]
      {The orthogonal subcategory problem in homotopy theory}
\author{Carles Casacuberta}
\author{Boris Chorny}
\thanks{During the preparation of this paper, the second-named author was a fellow of the Marie Curie 
Training Site with reference no.\ HPMT-CT-2000-00075 of the European Commission, hosted by the CRM (Barcelona).
The first-named author is currently supported by the Spanish Ministry of Education and Science
under MEC-FEDER grant MTM2004-03629.}

\address{Departament d'\`{A}lgebra i Geometria, Universitat de Barcelona, Gran Via de les Corts Catalanes, 585, \mbox{08007 Barcelona}, Spain}
\email{carles.casacuberta@ub.edu}

\address{Department of Mathematics, Middlesex College, The University of Western Ontario, London, Ontario N6A 5B7, Canada}
\email{bchorny2@uwo.ca}

\date{\today}
\dedicatory{} \commby{}

\begin{abstract}
It is known that the existence of localization with respect to an arbitrary (possibly proper) class of maps 
in the category of simplicial sets is implied by a large-cardinal axiom called \Vopenkas principle.
In this article we extend the validity of this result to any left proper,
combinatorial, simplicial model category $\cat M$
and show that, under additional assumptions on $\cat M$, every homotopy idempotent
functor is in fact a localization with respect to some set of maps. These results are valid
for the homotopy category of spectra, among other applications.
\end{abstract}

\maketitle
\section*{Introduction}

Localizing with respect to sets of maps is a common technique
in homotopy theory, as well as in other areas of Mathematics.
However, localizing with respect to proper classes of maps is a
more delicate issue, since the standard methods may fall into
set-theoretical difficulties. Due to these, for instance,
it is still unknown whether the existence of arbitrary cohomological
localizations can be proved or not using the ZFC axioms
(Zermelo--Fraenkel axioms with the axiom of choice).

An interesting step was made in \cite{CSS}, using results from
Ad\'amek and Rosick\'y \cite{AR}, by showing that a certain
large-cardinal axiom, called \Vopenkas principle, implies the 
existence of localization with respect to any (possibly proper)
class of maps in the category of simplicial sets. Hence,
\Vopenkas principle implies the existence of cohomological
localizations, among other consequences.

In this article we generalize the main results of \cite{CSS}
in two directions. In Section~1 we show that \Vopenkas principle implies 
the existence of localization with respect to any class of maps 
in left proper combinatorial simplicial model categories, and that any 
such localization is equivalent to localization with respect to a set
of maps. This also follows from results obtained by Rosick\'y and Tholen 
in \cite[\S 2]{RT} using a different argument.
The term ``combinatorial'', due to Jeff Smith, 
means locally presentable and
cofibrantly generated; see \cite{AR} and \cite{Hirschhorn} for the
definitions of these concepts.
Many useful model categories are combinatorial. Besides
the model category of simplicial sets, our results apply to the 
model category of symmetric spectra based on simplicial sets, 
and also to the model category of groupoids, among a number of
other cases.

Secondly, we address a closely related question, raised by Dror Farjoun in
\cite{Farjoun:v1}, asking if any functor $L$ that is idempotent
up to homotopy is equivalent to localization with
respect to some map $f$. He himself showed in \cite{Farjoun:JPAA}
that, if $L$ is a homotopy idempotent {\it simplicial} functor in the
category of simplicial sets, then it is equivalent to localization 
with respect to a {\it proper class} of maps. This result was improved
in \cite{CSS} by showing that the assumption that $L$ be simplicial
is unnecessary, and that, under \Vopenkas principle, the proper
class of maps defining $L$ can indeed be replaced by a set. Furthermore,
it was shown that such a replacement of a class by a set cannot be
done in general using only the ZFC axioms, since a counterexample was 
exhibited under a different axiom that is consistent with ZFC.

In Section~2 we show (without resorting to large-cardinal principles) 
that every homotopy idempotent
functor $L$ in a simplicial model category $\cat M$ 
is equivalent to localization with respect to a proper class of maps, assuming either 
that $L$ is simplicial or that $\cat M$ satisfies certain hypotheses stated in
\cite{RSS}, which allow to approximate any homotopy functor by
a simplicial functor. The hypotheses of \cite{RSS} hold in cofibrantly generated
stable model categories and in many other cases. Furthermore,
if one assumes that \Vopenkas principle is true and $\cat M$ is combinatorial,
then, again, the proper class of maps
defining $L$ can be replaced by a set. In most cases of interest, 
such a set of maps can be replaced by a single map (by taking the coproduct
of all maps in the set), but not always, as we show by means of an example
at the end of the paper.

In \cite{PhDII} an example was given of a homotopy idempotent functor in a
locally presentable (but not cofibrantly generated) model category which
fails to be a localization with respect to any set of maps. Hence our
results in Section~2 below are sharp.

\bigskip\noindent{\bf Acknowledgements:} \
Discussions with Mark Hovey, Jiri Rosick\'y, Brooke Shipley and Jeff Smith 
are greatly appreciated.

\section{Simplicial orthogonality}

Model categories were introduced by Quillen in \cite{Quillen} and have
recently been discussed in the books \cite{DHKS}, \cite{GJ}, \cite{Hirschhorn}, 
\cite{Hovey}, among many other places, with slight changes in the 
terminology and even in the assumptions. In this article we will assume
that model categories are complete, cocomplete, and equipped
with functorial factorizations. See \cite[\S~9]{DHKS},
\cite[\S~7]{Hirschhorn}, or \cite[\S~1]{Hovey} for full details.

Although our main results are stated for {\it simplicial} model categories
(for the definition, see for example \cite[II.3]{GJ} or \cite[9.1.5]{Hirschhorn}), several
of our steps require only {\it homotopy function complexes}, as introduced
in \cite{DK} and discussed in \cite[\S~5]{Hovey} or \cite[Ch.\ 17]{Hirschhorn}.
Thus, for any given model category $\cat M$, we make a
functorial choice of a fibrant simplicial set $\map(X,Y)$ for each $X$ and $Y$
in $\cat M$, whose homotopy type is the same as the diagonal of the bisimplicial set
${\cat M}(X^*,Y_*)$ where $X^*\to X$ is a cosimplicial resolution of $X$
and $Y\to Y_*$ is a simplicial resolution of $Y$. If $\cat M$ is a simplicial
model category, then any such choice of $\map(X,Y)$ is weakly equivalent to the 
simplicial set $\Map(X,Y)$ given by the simplicial structure of $\cat M$, provided that $X$ 
is cofibrant and $Y$ is fibrant. In fact, if $Q$ is a functorial cofibrant
approximation and $R$ is a functorial fibrant approximation in a simplicial
model category $\cat M$, then $\Map(QX,RY)$ is a good choice of a 
functorial homotopy function complex.

Before discussing simplicial orthogonality in model categories, we recall
the following older concepts from general category theory.
If $\cat C$ is any category, an object $X$ and a morphism $f\colon A\to B$
are called {\it orthogonal} (see \cite{AR} or \cite{CPP} for details
and motivation) if the induced function
\[
f^*\colon {\cat C}(B,X)\longrightarrow {\cat C}(A,X)
\]
is bijective. (We denote by ${\cat C}(X,Y)$ the set of morphisms from $X$ to $Y$
in $\cat C$.) If $L$ is an endofunctor of $\cat C$ equipped with a natural
transformation $\eta\colon\Id\to L$ such that $\eta L=L\eta$ and
$L\eta\colon L\to LL$ is an isomorphism, then $L$ is called an {\it idempotent functor}
or a {\it localization}. Then every object isomorphic to $LX$ for some $X$
is orthogonal to every map $f$ such that $Lf$ is an isomorphism, and in fact
these two classes determine
each other by the orthogonality relation; that is, an object is isomorphic
to $LX$ for some $X$ if and only if it is orthogonal to all maps $f$
such that $Lf$ is an isomorphism, and conversely.

As a special case, this terminology applies to the homotopy category $\HocatM$
associated with any model category $\cat M$.
Thus, orthogonality in $\HocatM$ between an
object $X$ and a map $f\colon A\to B$ amounts to the condition that
\begin{equation}
f^*\colon [B,X]\longrightarrow [A,X]
\end{equation}
be bijective, where $[X,Y]$ means, as usual, $\HocatM (X,Y)$.
Examples of idempotent functors in the homotopy category
of simplicial sets, such as homological localizations, have been studied 
since several decades ago; see, e.g., \cite{Bous:homology}.

Throughout the extensive study of localizations undertaken since then 
in homotopy theory, a stronger notion of ``simplicially enriched orthogonality'' 
came to be considered. There is no widely agreed terminology for it yet.
It was called {\it simplicial orthogonality} in \cite{CSS} and
{\it homotopy orthogonality} in \cite[\S 17]{Hirschhorn}. 
Thus, if
$\cat M$ is any model category with a functorial choice
of homotopy function complexes, an object $X$ and a map $f\colon A\to B$
will be called {\it homotopy orthogonal} or {\it simplicially orthogonal} 
if the induced map of simplicial sets
\begin{equation}
f^*\colon \map(B,X)\longrightarrow \map(A,X)
\end{equation}
is a weak equivalence. 
Since there is a natural bijection between $\pi_0\map(X,Y)$ and $[X,Y]$,
homotopy orthogonality implies orthogonality in $\HocatM$ in the sense of~(1).
Although plenty of examples show that the converse is not true,
we discuss in Section~2 an important situation where the converse holds.

The fibrant objects $X$ that are homotopy orthogonal to a given map $f$
are usually called {\it $f$-local}. More generally, if $\calS$ is any
class of maps, we denote by ${\calS}^{\hbot}$ the class
of fibrant objects that are homotopy orthogonal to all the maps in $\calS$,
and say that objects in ${\calS}^{\hbot}$ are {\it $\calS$-local}.
The homotopy orthogonal complement of a class $\calD$ of~objects
is defined similarly.

For a class of morphisms $\calS$ in a model category $\cat M$,
an {\it $\calS$-localization} is an endofunctor 
$L\colon \cat M \rarrow \cat M$ equipped with a natural transformation 
$\eta\colon \Id \rarrow L$ such that $\eta L \simeq L\eta$ and $L\eta\colon L\to LL$ 
is a weak equivalence on all objects, $\eta_X\colon X\to LX$ is in $({\calS}^{\hbot})^{\hbot}$
and $LX$ is $\calS$-local for all $X$.
We also call it a {\it localization with respect to $\calS$}
or, generically, a {\it homotopy localization}.
Indeed, every homotopy localization yields a localization functor in $\HocatM$.

The {\it orthogonal subcategory problem} in category theory 
asks if, given a class $\calS$ of morphisms in a
category $\cat C$, the class of objects orthogonal to $\calS$ 
is the class of local objects of some localization functor.
Sufficient conditions for an affirmative answer, depending of course
on the properties of the category $\cat C$, can be found in
many places; see for instance \cite{AR}. The extent to which
the orthogonal subcategory problem in locally presentable categories
depends on set-theoretical axioms was previously discussed in
\cite{ART}.

The homotopy version of the orthogonal
subcategory problem, that we address in this article, 
consists of asking the same question in a model 
category $\cat M$, by considering homotopy orthogonality as given by~(2). 
The reason for using (2) instead of (1) is the
fact that the orthogonal subcategory problem in $\HocatM$, using (1)
as orthogonality relation, very often has a negative answer.
For instance, there is no localization in the homotopy category of
simplicial sets onto the class of simply connected spaces. 
See \cite{CC:PAMS} for a more elaborate counterexample.

It is well known that the homotopy orthogonal subcategory problem
has a solution if $\calS$ is a set and $\cat M$ satisfies certain assumptions, which 
vary slightly depending on the authors.
We will call a model category {\it combinatorial}
if it is locally presentable and cofibrantly generated.
The definition of a locally presentable category can be found
in \cite{AR}, and the definition of a cofibrantly generated model category
is contained, e.g., in \cite{Hirschhorn}. The notion of properness is also discussed
in \cite{Hirschhorn}.

\begin{theorem}\label{localization}
Let $\cat M$ be a left proper, combinatorial, simplicial model category.
For any set of maps $\calS$ there is a localization functor $L$
with respect to $\calS$.
\end{theorem}
\begin{proof}
The core of the proof is in \cite{Bous:factor}. See \cite{Hirschhorn} for an updated approach.
\end{proof}

As far as we know, there is no way to prove this when $\calS$ is a proper class,
not even for simplicial sets, using the ordinary axioms of set theory.
In \cite{CSS} it was shown that the statement of Theorem~\ref{localization} holds
for a proper class $\calS$ in the model category of simplicial sets 
using a suitable large-cardinal axiom (\Vopenkas principle).
We now undertake a generalization of this fact to other model categories.

Recall that a partially ordered set $A$ is called {\it $\lambda$-directed},
where $\lambda$ is a regular cardinal, if every subset of $A$ of cardinality
smaller than $\lambda$ has an upper bound.

\begin{lemma}\label{cocomplete}
Let $\calD$ be any class of objects in a combinatorial simplicial model category $\cat M$,
and let ${\calS}$ be its homotopy orthogonal complement.
Then there exists a regular cardinal $\lambda$ such that $\calS$ is closed under 
$\lambda$-directed colimits in the category of maps of $\cat M$.
\end{lemma}

\begin{proof}
Let $I$ be a set of generating cofibrations for the model category $\cat M$. 
Choose a regular cardinal $\lambda$ such that any object of the set of domains and codomains of maps in $I$ 
is $\lambda$-presentable (such a cardinal exists since the category $\cat M$ is locally presentable). 
Let $A$ be any $\lambda$-directed partially ordered set, and suppose given
a diagram $f\colon A\to\cat \ArrM$, where $\ArrM$ is the category of maps in $\cat M$.
Let us depict it, for simplicity, as a chain:
\begin{equation}
\begin{CD}
 X_0       @>>>    X_1  @>>>  \cdots     @>>>    X_n    @>>> \cdots  \\
@V{f_0}VV        @V{f_1}VV       @.          @V{f_n}VV          @.   \\
 Y_0       @>>>    Y_1  @>>>  \cdots     @>>>    Y_n    @>>> \cdots.  \\
\end{CD}
\end{equation}
Suppose that the maps $f_i$ are in $\calS$ for each $i\in A$. Since $\cat M$ is cocomplete,
$\ArrM$ is cocomplete as well, and we may consider the colimit of the diagram $f$.
We need to show that the induced map
$\colim f_i \colon \colim X_i \longrightarrow \colim Y_i$
is also in $\calS$.

Consider the category $\cat M^A$ of $A$-indexed diagrams in $\cat M$, and endow it with a 
model structure as described in \cite[11.6]{Hirschhorn}. Thus, weak equivalences 
and fibrations are objectwise, 
and cofibrations are retracts of free cell complexes. The diagram (3) may be viewed as a single map 
in $\cat M^A$. Apply the cofibrant approximation functor to this map using the above model structure,
hence obtaining the following commutative diagram in $\cat M$:
\[
\def\labelstyle{\scriptstyle}
\xymatrix@!0{
  \tilde X_0  \ar[dd]_<(.6){\tilde f_0} \ar@{->>}[dr]^{\dir{~}} \ar@{^{(}->}[rr]  &      &  \tilde X_1 \ar'[d][dd]_<(.3){\tilde f_1} \ar@{->>}[dr]^{\dir{~}} \ar@{^{(}->}[rr]        &   & \cdots \ar@{^{(}->}[rr]&  & \tilde X_n \ar'[d][dd]_{\tilde f_n} \ar@{->>}[dr]^{\dir{~}} \ar@{^{(}->}[rr] &  & \cdots\\
         &  X_0 \ar[dd]^<(.6){f_0} \ar[rr]               &      &  X_1  \ar[dd]^<(.6){f_1} \ar[rr]        &   & \cdots \ar[rr]&  & X_n \ar[rr]\ar[dd]^<(.6){f_n} &  & \cdots\\
   \tilde Y_0  \ar@{->>}[dr]^{\dir{~}}\ar@{^{(}->}'[r][rr]      &      &  \tilde Y_1 \ar@{->>}[dr]^{\dir{~}} \ar@{^{(}->}'[r][rr]       &   & \cdots \ar@{^{(}->}[rr]&  & \tilde Y_n \ar@{->>}[dr]^{\dir{~}} \ar@{^{(}->}'[r][rr] &  & \cdots\\
         &  Y_0 \ar[rr]                       &      &  Y_1  \ar[rr]      &   & \cdots \ar[rr]        &  &   Y_n \ar[rr] &  & \cdots, \\
}
\]
where $\tilde f_i$ is a cofibrant approximation of $f_i$. 

For every $Z\in {\calD}$, let $\hat Z$ be a fibrant approximation to $Z$. The induced map
\[
\Map(\colim \tilde f_i, \hat Z)\colon \Map(\colim \tilde Y_i,
\hat Z) \longrightarrow \Map(\colim \tilde X_i, \hat Z)
\]
can be written as a limit of maps
\[
\lim \Map( \tilde f_i, \hat Z)\colon \lim \Map(\tilde Y_i,
\hat Z) \longrightarrow \lim \Map(\tilde X_i, \hat Z),
\]
each of which is a fibration of simplicial sets, since $\tilde f_i$ is a cofibration
for all $i$. Hence, $\Map(\colim\tilde f_i,\hat Z)$ is a homotopy inverse limit of weak equivalences,
so it is itself a weak equivalence. This shows that $\colim\tilde f_i$ is in $\calS$.

Trivial fibrations in $\cat M$ are preserved 
under $\lambda$-directed colimits, since the set of generating cofibrations has $\lambda$-presentable 
domains and codomains. From the commutative diagram
\[
\xymatrix{
 \colim \tilde X_i \ar@{->>}[r]^\sim \ar[d]_{\colim \tilde f_i}      
 &  \colim X_i\ar[d]^{\colim f_i} \\
 \colim \tilde Y_i \ar@{->>}[r]^\sim  &  \colim Y_i\\
 }
\]
we conclude that the map $\colim \tilde f_i$ is a cofibrant approximation of the
map $\colim f_i$. Hence, $\colim f_i$ is in $\calS$, as claimed.
\end{proof}

The statement of \Vopenkas principle and enough motivation for its use in this context can
be found in \cite{AR}, \cite{CSS}, and \cite{RT}.

\begin{lemma}\label{set-of-generators}
Suppose that \Vopenkas principle is true. Let $\calD$ be any class of objects in a 
combinatorial simplicial model category $\cat M$, and let ${\calS}={\calD}^{\hbot}$. 
Then there exists a set of maps $\calX$ such that ${\calX}^{\hbot} = {\calS}^{\hbot}$.
\end{lemma}

\begin{proof}
By abuse of notation, we also denote by $\calS$ the full subcategory of $\ArrM$ 
generated by the class $\calS$. Since $\cat M$ is locally presentable, 
$\ArrM$ is also locally presentable. Then, assuming \Vopenkas principle, it follows from
\cite[Theorem~6.6]{AR} that $\calS$ is bounded, i.e., it has a small dense subcategory. 
We have shown in Lemma~\ref{cocomplete} that there exists a regular cardinal $\lambda$ such that $\calS$ is 
closed under $\lambda$-directed colimits in the category $\ArrM$. Hence, by \cite[Corollary~6.18]{AR}, 
the full subcategory generated by $\calS$ in $\ArrM$ is accessible. Thus, for a certain regular 
cardinal $\lambda_0 \geq \lambda$, the class $\calS$ contains a set $\calX$ of $\lambda_0$-presentable objects
such that every object of $\calS$ is a $\lambda_0$-directed colimit of objects of $\calX$.

Since $\calX \subset \calS$, we have  $\calX^\hbot \supset \calS^\hbot$ and $(\calX^\hbot)^\hbot 
\subset (\calS^\hbot)^\hbot = \calS$. Our aim now is to show the reverse inclusion 
$(\calX^\hbot)^\hbot \supset \calS$. By Lemma~\ref{cocomplete}, $(\calX^\hbot)^\hbot$ 
is closed under $\lambda$-directed colimits .
Hence $(\calX^\hbot)^\hbot$ is also closed under $\lambda_0$-directed colimits and  
every element of $\calS$ is a $\lambda_0$-directed colimit of elements of $\calX$. Then we can choose 
$\calX$ as our generating set.
\end{proof}

\begin{theorem}\label{HOSP}
Let $\cat M$ be a left proper, combinatorial, simplicial model category.
If \Vopenkas principle is assumed true, then
for any (possibly proper) class of maps $\calS$ there is a localization 
functor $L$ with respect to $\calS$.
\end{theorem}
\begin{proof}
By Lemma~\ref{set-of-generators}, there exists a set $\calX$ of maps in $\cat M$ such that 
${\calX}^{\hbot}={\calS}^{\hbot}$. 
Then the localization with respect to $\calX$, which exists by Theorem~\ref{localization},
has ${\calS}^{\hbot}$ as its class of local objects.
\end{proof}

The statement of Theorem~\ref{HOSP} is our positive solution of the 
homotopy orthogonal subcategory problem
in suitable model categories, under \Vopenkas principle.

\section{Idempotent functors and simplicial orthogonality}\label{Farjoun-problem}

The next theorem is motivated by results of Dror Farjoun in \cite{Farjoun:JPAA},
and generalizes the argument given for groupoids in \cite[Theorem 2.4]{CGT}.
As in the previous section, we assume that a functorial choice of a homotopy
function complex has been made in every model category under consideration.

\begin{theorem}\label{continuity}
Let $\cat M$ be any model category.
Let $L$ be an endofunctor in the homotopy category $\HocatM$ with the following properties:
\begin{itemize}
\item[(a)] There is a natural transformation $\eta\colon\Id\to L$ in $\HocatM$
such that $L\eta\simeq\eta L$ and $L\eta\colon L\to LL$ is a weak equivalence on all objects.
\item[(b)] There is a map $l_{X,Y}\colon \map(X,Y)\to \map(LX,LY)$ for all $X$, $Y$,
which is natural in both variables in $\HocatM$.
\item[(c)] $(\eta_X)^*\circ l_{X,Y}\simeq (\eta_Y)_*$ for all $X$ and $Y$.
\end{itemize}
Then $(\eta_X)^*\colon \map(LX,LY)\to \map(X,LY)$
is a weak equivalence for all $X$, $Y$.
\end{theorem}

\begin{proof}
Let us write $Z=LY$ for simplicity.
The assumption (a) says precisely that $L$ is idempotent in the
homotopy category $\HocatM$. Hence, among other consequences of this fact, 
$\eta_Z\colon Z\to LZ$ is a weak equivalence.
Let $\xi\colon \map(LX,LZ)\to \map(LX,Z)$ be a homotopy inverse of $(\eta_{Z})_*$.
We claim that $\xi\circ l_{X,Z}$ is a homotopy inverse of $(\eta_X)^*$. On one hand,
\[
\xi\circ l_{X,Z}\circ (\eta_X)^*\simeq \xi\circ (L\eta_X)^*\circ l_{LX,Z}
\simeq \xi\circ (\eta_{LX})^*\circ l_{LX,Z}\simeq \xi\circ (\eta_Z)_*\simeq\identity.
\]
(The first homotopy equivalence comes from the naturality of $l$, and the third 
homotopy equivalence is given by (c).)
On the other hand,
$(\eta_X)^*\circ\xi \circ l_{X,Z} \simeq \xi\circ (\eta_Z)_*\circ
(\eta_X)^*\circ\xi\circ l_{X,Z} \simeq
\xi\circ (\eta_X)^*\circ (\eta_Z)_*\circ\xi\circ l_{X,Z} \simeq
\xi\circ (\eta_X)^*\circ l_{X,Z} \simeq \xi\circ (\eta_Z)_* \simeq
\identity$, which completes the proof.
\end{proof}

This result says that, under the stated hypotheses,
$L$ is localization with respect to the class of maps $\eta_X$ for all $X$.
Indeed, the homotopy orthogonal complement of this class of maps contains
all objects of the form $LY$ by Theorem~\ref{continuity}, and if $Z$ is
any object such that $\map(LX,Z)\simeq\map(X,Z)$ for all $X$, then
$[LX,Z]\cong[X,Z]$ for all $X$ as well, from which it follows that $Z\simeq
LY$ for some $Y$, since $L$ is an idempotent functor in $\HocatM$.

Assumptions (b) and (c) in Theorem~\ref{continuity} need not be satisfied
by arbitrary idempotent functors in $\HocatM$, not even by
those derived from functors in $\cat M$.
A stronger form of assumption (b) played a central role in Dror Farjoun's approach 
in \cite{Farjoun:JPAA}. A~way to avoid this continuity hypothesis in the model 
category of simplicial sets was found in \cite{CSS}, by using a convenient
homotopy function complex.

Here we use another method, based on Proposition~6.4 in \cite{RSS}.
If $\cat M$ is a model category, let $s\cat M$ denote the category of
simplicial objects over $\cat M$. The {\it canonical model structure}
on $s\cat M$ is the only one where every level equivalence is a weak equivalence,
the cofibrations are the Reedy cofibrations, and the fibrant objects are the homotopically
constant Reedy fibrant objects (see \cite{RSS}
for motivation and further details). This model structure need not exist; however,
when it exists, $s\cat M$ is a simplicial model category that is Quillen
equivalent to $\cat M$. Sufficient conditions for its existence were given
in \cite{RSS}. 

Pointed model categories where the suspension functor and
the loop functor are inverse equivalences on the homotopy category are
called {\it stable}. According to \cite[Proposition 4.5]{RSS}, if $\cat M$ is
a proper, cofibrantly generated, stable model category, then the canonical
model structure on $s\cat M$ exists, and hence $\cat M$ is Quillen equivalent
to a simplicial model category, namely $s\cat M$.

A functor $F$ in a simplicial model category is called {\it simplicial}
if it comes equipped with natural maps of simplicial sets
$\Map(X,Y)\to\Map(FX,FY)$ preserving composition and identity;
see \cite[IX.1]{GJ} or \cite[9.8]{Hirschhorn}.

\begin{theorem}
Let $\cat M$ be a cofibrantly generated simplicial model category where the
canonical model category structure exists in $s\cat M$.
Let $L$ be an endofunctor in $\cat M$ equipped with a natural transformation
$\eta\colon\Id\to L$. Suppose that $L$ preserves weak
equivalences and is idempotent
in the homotopy category $\HocatM$. Then $L$ is equivalent to homotopy
localization with respect to some class of maps.
\end{theorem}

\begin{proof}
Let $F$ denote a functorial simplicial approximation, as given
by \cite[Proposition 6.4]{RSS}. Thus, $(FL)X=|Q\hat L\,\Sing X|$
for each object $X$, where $(\Sing X)_n=X^{\Delta[n]}$ for all $n$, 
the realization functor $|-|$
is its left adjoint, $\hat L$ is the dimensionwise prolongation of $L$ over
$s\cat M$, and $Q$ is a simplicial cofibrant replacement functor.
By its construction, $FL$ is a simplicial functor
in $\cat M$ which comes with a ziz-zag of weak equivalences between $L$ and $FL$.

Then there is also a natural transformation of simplicial functors
$F\eta\colon F\Id\to FL$, where $F\Id$ need not be the identity.
As detailed in \cite[IX.1]{GJ}, this implies that there are natural maps
\[
{\alpha}_{X,Y}\colon\Map(X,Y)\longrightarrow\Map((F\Id)X,(F\Id)Y)
\]
and
\[
{\beta}_{X,Y}\colon\Map(X,Y)\longrightarrow\Map((FL)X,(FL)Y)
\]
such that 
\begin{equation}
((F\eta)_X)^* \circ {\beta}_{X,Y} = ((F\eta)_Y)_* \circ {\alpha}_{X,Y}
\end{equation}
for all $X$ and $Y$.

Let us check that the derived functor $FL$ satisfies the assumptions
of Theorem~\ref{continuity} in the homotopy category $\HocatM$.
Since $FL\simeq L$ and $F\Id\simeq\Id$, there is indeed a natural 
transformation $\Id\to FL$ in $\HocatM$ rendering $FL$ into an
idempotent functor, hence fulfilling assumption~(a). 
Next, since ${\alpha}_{X,Y}$ is a weak equivalence for all $X$ 
and~$Y$, we may choose a homotopy inverse and consider
$l_{X,Y}={\beta}_{X,Y}\circ{\alpha}_{X,Y}^{-1}$. This map gives rise to
a corresponding map of homotopy function complexes,
as required in~(b). Finally, the assumption (c) follows from~(4).

Then the conclusion of Theorem~\ref{continuity}
implies that $FL$, and hence also $L$, is equivalent to localization 
with respect to the class of maps of the form $\eta_X$ for all $X$.
\end{proof}

Now the results of the previous section yield the following
answer to Dror Farjoun's problem in sufficiently good model categories.

\begin{theorem}\label{Farjoun-theorem}
Assuming \Vopenkas principle, any homotopy idempotent functor in a
combinatorial simplicial model category $\cat M$ where the canonical
model structure exists in $s\cat M$
is equivalent to localization with respect to some set of maps.
\end{theorem}

This result applies, as an important case, to the
stable homotopy category of Adams--Boardman (using, for example,
the model category of symmetric spectra based on simplicial sets).

In the category of simplicial sets, the set $\calX$ of maps given by Theorem~\ref{Farjoun-theorem}
can be replaced by a single map $f$, namely the coproduct $\coprod_{g\in\calX} g$ of all maps in~$\calX$. 
In a general model category one has to be more careful, in view of the next counterexample.

Consider the model category which is a product of two copies of the category of
spaces, i.e., the category of diagrams of spaces over the discrete category
with two objects equipped with the Bousfield--Kan model structure. 
Take $S = \{f,g\}$ for
\[
f\colon (\emptyset,\emptyset)\longrightarrow (\ast,\emptyset)\quad
\text{and}\quad g\colon (\emptyset,\ast)\longrightarrow
\left(\emptyset,\ast{\textstyle\coprod} \ast\right).
\]
An object $(X,Y)$ is $S$-local if and only if $X$ and $Y$ are fibrant, $X$
is contractible and $Y$ is either contractible or empty.

Suppose that there exists a map
\[
h\colon (A,B)\longrightarrow (C,D)
\]
such that any $S$-local object is also $h$-local, and vice versa. The object
$(X,\emptyset)$ is $h$-local if and only if $X$ is contractible. This
condition implies that both spaces $B$ and $D$ are empty; otherwise, for any
space $Z$, either contractible or not, the object $(Z,\emptyset)$ would be
$h$-local. But in this case any object $(X,Y)$ with contractible $X$
becomes $h$-local, hence the contradiction.


\end{document}